\documentclass{amsart}
\usepackage{amssymb}
\usepackage{amsmath}
\usepackage{latexsym}
\usepackage{eepic}
\usepackage{epsfig}
\usepackage{graphicx}
\usepackage{pb-diagram,lamsarrow,pb-lams,amscd}
\usepackage{eufrak}
\usepackage{mathrsfs}
\usepackage[dutch,english]{babel}

\DeclareMathAlphabet{\mathpzc}{OT1}{pzc}{m}{it}
\newtheorem{theorem}{Theorem}[section]
\newtheorem{theorem-definition}[theorem]{Theorem-Definition}
\newtheorem{lemma-definition}[theorem]{Lemma-Definition}
\newtheorem{definition-prop}[theorem]{Proposition-Definition}

\newtheorem{prop}[theorem]{Proposition}
\newtheorem{lemma}[theorem]{Lemma}

\newtheorem{definition}[theorem]{Definition}

\newtheorem{conjecture}[theorem]{Conjecture}

\newenvironment{remark}{\vspace{4pt}\noindent\textbf{Remark.}}{\qed\vspace{4pt}}

\newcommand{\Z}{\ensuremath{\mathbb{Z}}}
\newcommand{\Q}{\ensuremath{\mathbb{Q}}}

\newcommand{\C}{\ensuremath{\mathbb{C}}}

\newcommand{\X}{\ensuremath{\mathscr{X}}}

\newcommand{\mY}{\ensuremath{\mathfrak{Y}}}

\renewcommand{\C}{\ensuremath{\mathbb{C}}}

\renewcommand{\X}{\ensuremath{\mathfrak{X}}}

\renewcommand{\mY}{\ensuremath{\mathfrak{Y}}}
\newcommand{\mZ}{\ensuremath{\mathfrak{Z}}}

\newcommand{\Spec}{\ensuremath{\mathrm{Spec}\,}}
\newcommand{\Spf}{\ensuremath{\mathrm{Spf}\,}}

\numberwithin{equation}{section}

\begin{document}
\title{Gauss-Manin connection and $t$-adic geometry}
\author[Johannes Nicaise]{Johannes Nicaise}
\address{Universit\'e Lille 1\\
Laboratoire Painlev\'e, CNRS - UMR 8524\\ Cit\'e Scientifique\\59655 Villeneuve d'Ascq C\'edex \\
France} \email{johannes.nicaise@math.univ-lille1.fr} \maketitle
\section{Introduction}
Let $k$ be a field of characteristic zero, denote by $R$ the ring $k[[t]]$ of formal power series over $k$, and by $K$
the field $k((t))$ of Laurent series
over $k$. The aim of this note is to show that the de Rham cohomology of any separated and smooth rigid $K$-variety $X$
carries a natural formal meromorphic connection $\partial_{X}$, which we call the Gauss-Manin connection.

It relates to the classical Gauss-Manin connection in the following way. Let
 $S$ be a smooth $k$-curve, $0$ a point of $S(k)$, and $t$ a local parameter on $S$ at $0$.
 We put $S^o=S- \{0\}$.
 The choice of $t$ defines
 a morphism of $k$-schemes $\widehat{\eta}:\Spec K\rightarrow S^o$. Let $f:Y\rightarrow S^o$
be a proper and smooth morphism, and put $X=Y\times_{S^o} \widehat{\eta}$. If we denote by $\nabla$
the classical Gauss-Manin connection
on $R^if_*(\Omega^{\bullet}_{Y/S^o})$, then the covariant derivative
$\partial_t=\nabla_{\!\!\frac{\partial}{\partial t}}$ induces a formal
meromorphic connection on the de Rham cohomology space $H_{dR}^i(X/K)$ for each $i\geq 0$.
We showed that the natural GAGA isomorphism $H_{dR}^i(X/K)\cong H_{dR}^i(X^{an}/K)$
commutes with the connections $\partial_t$ and $\partial_{X^{an}}$ (Theorem \ref{compar-prop}).

In the local case, we conjecture that a similar comparison result holds. Let $f:\C^{n+1}\rightarrow \C$ be
a complex analytic map with an isolated singularity at $x\in f^{-1}(0)$ and denote by $\mathscr{F}_x$ the analytic Milnor fiber
of $f$ at $x$. Our conjecture compares
the formal meromorphic connection $(H^*_{dR}(\mathscr{F}_x),\partial_{\mathscr{F}_x})$
to the Gauss-Manin connection on
the relative de Rham cohomology of the Milnor fibration of $f$ at $x$. See
Conjecture \ref{conj1} for a precise statement.

This note merely serves as an announcement of the principal results.
Detailed proofs will appear in a forthcoming
paper.
\subsection*{Notation}
Throughout this article, we denote by $k$ a field of characteristic zero, and we put
 $R=k[[t]]$ and $K=k((t))$. We endow $R$ with its $t$-adic topology, and we fix a $t$-adic absolute value $|\cdot|$
 on $K$ (determined by the choice of $|t|\in \,]0,1[$\,).

To simplify arguments, all rigid $K$-varieties will be assumed to be separated. If $X$ is a separated
 smooth rigid $K$-variety,
 we denote by $H_{dR}^*(X/K)$ its ``na\"ive'' de Rham cohomology, i.e. the hypercohomology of the de Rham complex
 $\Omega^\bullet_{X/K}$. In the cases of most importance to us, $X$ will be partially proper so that na\"ive de
 Rham cohomology coincides with overconvergent de Rham cohomology \cite[1.8(b)]{rigiddeRham}. If $X$ is a smooth algebraic $K$-variety,
 we also denote its de Rham cohomology by $H_{dR}^*(X/K)$.

 Let $\X$ be a separated  Noetherian adic formal scheme, endowed with a morphism $\X\rightarrow \Spf R$. We denote
 by $\X_0$ the reduction of $\X$, i.e. the closed subscheme defined by the largest ideal of definition. This is a separated
 reduced $k$-scheme.
 We say that $\X$ is a special formal $R$-scheme if $\X_0$ is of finite type over $k$, and that
$\X$ is $stft$ if $\X$ is topologically of finite
 type over $R$. Any $stft$ formal $R$-scheme is special.

We denote by $(SpF/R)$ the category of special formal $R$-schemes, by
$$(\cdot)_\eta:(SpF/R)\rightarrow (Rig/K):\X\mapsto \X_\eta$$ the generic fiber functor to
the category $(Rig/K)$ of separated
rigid $K$-varieties, and by $sp_{\X}:\X_\eta\rightarrow \X$ the natural morphism of ringed sites \cite[0.2.6]{bert}.
We say that $\X$ is generically smooth if $\X_\eta$ is a smooth rigid $K$-variety.

 For any special formal $R$-scheme $\X$, we
  consider the natural functor
 $$(\cdot)_{rig}=(sp_{\X})^*:(Mod_{\X})\rightarrow (Mod_{\X_\eta})$$ from the category of
 $\mathcal{O}_{\X}$-modules to the category
 of $\mathcal{O}_{\X_\eta}$-modules. This functor is exact and takes coherent $\mathcal{O}_{\X}$-modules to coherent
 $\mathcal{O}_{\X_\eta}$-modules.
 It is left adjoint to the direct image functor
$$(sp_\X)_*:(Mod_{\X_\eta})\rightarrow (Mod_{\X})$$
Since $\X_\eta$ is quasi-Stein if $\X$ is affine \cite[2.4]{Ni-trace}, Kiehl's Theorem B \cite[2.4]{Kiehl} implies that
the restriction of $(sp_{\X})_*$ to the category of coherent $\mathcal{O}_{\X_\eta}$-modules is exact for any special formal
$R$-scheme $\X$.
  If $h:\mY\rightarrow \X$ is a morphism of special formal $R$-schemes and $M$ is a $\mathcal{O}_{\X}$-module,
 there is a canonical isomorphism $(h^*M)_{rig}\cong (h_\eta)^*M_{rig}$ since $h\circ sp_{\mY}=sp_{\X}\circ h_\eta$.

 We denote by
 $$(\cdot)^h:(sft/\C)\rightarrow (An/\C):X\mapsto X^h$$
 the GAGA functor from the category of separated $\C$-schemes of finite type to the category of complex
 analytic spaces, and by
 $$(\cdot)^{an}:(sft/K)\rightarrow (Rig/K):X\mapsto X^{an}$$ the GAGA functor from the category of separated $K$-schemes
 of finite type to the category of separated rigid $K$-varieties.

If $A$ is any ring and $M^\bullet$ is a complex of $A$-modules, then we denote by
$\mathcal{H}^i(M^\bullet)$ the $i$-th cohomology space of $M^\bullet$, for each $i\in \Z$.
\section{The Gauss-Manin connection of a smooth
rigid $K$-variety}\label{sec-gm} \subsection{Differentials over the residue
field} Let $X$ be a separated smooth rigid
variety over $K$. We will define, for each $i\geq 0$, a coherent
$\mathcal{O}_X$-module $\Omega^i_{X/k}$ of differential $i$-forms
over $k$, and a de Rham complex $(\Omega^\bullet_{X/k},d)$.

First, consider the case where $X$ is quasi-compact, and choose a
$stft$ formal $R$-model $\X$ for $X$, i.e. a $stft$ formal $R$-scheme $\X$ endowed with an isomorphism
 $\X_\eta\cong X$. Since
$\X$ is of pseudo-finite type over $k$ in the terminology of \cite{formal1}, we have coherent
$\mathcal{O}_{\X}$-modules $\Omega^i_{\X/k}$ of continuous differential
$i$-forms over $k$, which fit into a de Rham complex
$(\Omega^\bullet_{\X/k},d)$. We put
$\Omega^i_{X/k}=(\Omega^i_{\X/k})_{rig}$ for each $i\geq 0$. This
definition is independent of the chosen model, by the fact that any pair of formal models can be
dominated by a third \cite[4.1]{formrigI}, and by the following
lemma.

\begin{lemma}\label{blup}
If $h:\mZ\rightarrow \mY$ is a morphism of generically smooth special formal
$R$-schemes, and $h_\eta$ is \'etale, then the natural map
$$h_\eta^*((\Omega^i_{\mY/k})_{rig})\cong(h^*\Omega^i_{\mY/k})_{rig}\rightarrow
(\Omega^i_{\mZ/k})_{rig}$$ is an isomorphism.
\end{lemma}
\begin{proof}
This follows from \cite[7.18-19]{Ni-trace}.
\end{proof}

Combined with \cite[4.1]{formrigI}, Lemma \ref{blup} also implies that for each $i\geq 0$ and any quasi-compact open subvariety $U$ of $X$,
 the restriction of
$\Omega^i_{X/k}$ to $U$ is canonically isomorphic to
$\Omega^i_{U/k}$. Note that for each $i\geq 0$
there is a canonical isomorphism  \begin{equation}\label{wedge}\Omega^i_{X/k}\cong
\bigwedge^i_{\mathcal{O}_X}\Omega^1_{X/k}\end{equation}

Now we define differentials $d:\Omega^i_{X/k}\rightarrow
\Omega^{i+1}_{X/k}$. Let $\X$ be a $stft$ formal $R$-model of $X$ and consider the sheaves
$$(sp_\X)_*\Omega^i_{X/k}\cong \Omega^i_{\X/k}\otimes_R K$$ for $i\geq 0$.
Applying the Leibniz rule one sees that there exists a unique
$k$-derivation $d_{\X}:(sp_\X)_*\mathcal{O}_{X}\rightarrow (sp_\X)_*\Omega^1_{X/k}$
 such that
the natural diagram
$$\begin{CD}
\mathcal{O}_{\X}@>>> (sp_{\X})_*\mathcal{O}_X
\\ @VdVV @VVd_{\X}V
\\ \Omega^1_{\X/k}@>>> (sp_{\X})_*\Omega^1_{X/k}\end{CD}$$
commutes.
 If $h:\mY\rightarrow \X$ is a morphism of $stft$ formal $R$-schemes such that $h_\eta$ is an isomorphism, then it
is easily seen that the natural isomorphisms $(sp_\X)_*\Omega^i_{X/k}\cong h_*(sp_\mY)_*\Omega^i_{X/k}$ $(i=1,2)$ commute
with the differentials $d_{\X}$ and $h_*d_{\mY}$. By \cite[4.4]{formrigI} and the fact that the categories of sheaves
w.r.t. the weak, resp.
strong $G$-topology on $X$ are equivalent, we can conclude that there exists
a unique map of sheaves $d:\mathcal{O}_{X/k}\rightarrow \Omega^{1}_{X/k}$ such that for any $stft$ formal
$R$-model $\X$ of $X$ we have $(sp_{\X})_*(d)=d_{\X}$. The map $d$ is again a $k$-derivation.

Using the natural isomorphisms in (\ref{wedge}), we obtain a de Rham complex
$(\Omega^\bullet_{X/k},d)$ for each quasi-compact separated smooth rigid $K$-variety $X$.
 If $X$ is a separated smooth rigid $K$-variety, we take an admissible cover of $X$
 by quasi-compact open subsets. If $W\subset V$ are quasi-compact open subsets of $X$, then
it is clear from our constructions that the de Rham complex
$(\Omega^\bullet_{W/k},d)$ is the restriction to $W$ of the de
Rham complex $(\Omega^\bullet_{V/k},d)$. Therefore, we can glue
our local constructions to coherent $\mathcal{O}_{X}$-modules $\Omega^i_{X/k}$ which fit
into a de Rham complex $(\Omega^\bullet_{X/k},d)$.

\begin{lemma}\label{special}
If $\X$ is a generically smooth special formal $R$-scheme, then there is a canonical
isomorphism
$$\Omega^i_{\X_\eta/k}\cong (\Omega^i_{\X/k})_{rig}$$ for each
$i\geq 0$, and the natural maps
$$\Omega^i_{\X/k}\rightarrow (sp_{\X})_*\Omega^i_{\X_\eta/k}$$ commute with the differentials
$d$ and $(sp_{\X})_*d$.
\end{lemma}
\begin{proof}
Let $\mathcal{I}$ be the largest ideal of definition on $\X$.
For each $n> 0$, we denote by $\pi(n):\X(n)\rightarrow \X$ the dilatation with center $(\mathcal{I}^n,t)$
(see \cite[2.20]{Ni-trace}). Then $\X(n)$ is a $stft$ formal $R$-scheme, the morphism
$\pi(n)_\eta:\X(n)_\eta\rightarrow \X_\eta$ is an open immersion,
 and the morphisms $\pi(n)_\eta$ form an admissible
cover of $\X_\eta$ by \cite[2.25]{Ni-trace}.
Moreover, for any coherent $\mathcal{O}_{\X}$-module $M$, the restriction of $M_{rig}$ to $\X(n)_\eta$
is canonically
isomorphic to $(\pi(n)^*M)_{rig}$. Now the statement follows immediately from Lemma \ref{blup} and
the definition of the complex
$(\Omega^\bullet_{\X_\eta/k},d)$.
\end{proof}
\begin{lemma}\label{etale}
If $h:Y\rightarrow X$ is an \'etale morphism of separated
smooth rigid $K$-varieties, then there exist natural isomorphisms
$$h^*\Omega^i_{X/k}\rightarrow \Omega^i_{Y/k}$$
for $i\geq 0$.
\end{lemma}
\begin{proof}
We may assume that $Y$ and $X$ are quasi-compact; then the
statement follows immediately from Lemma \ref{blup}.
\end{proof}

Let $X$ be a separated smooth rigid $K$-variety, and consider the section $dt$ in
$\Omega^1_{X/k}(X)$.
\begin{lemma}\label{exdef}
The wedge product $dt\wedge(\cdot )$ induces an exact sequence of
complexes
$$\begin{CD}
\Omega^\bullet_{X/k}[-1]@>dt\wedge(\cdot )>>\Omega^\bullet_{X/k}@>>>
\Omega^\bullet_{X/K}@>>> 0
\end{CD}$$
\end{lemma}
\begin{proof}
We may assume that $X$ has a $stft$ formal $R$-model $\X$.
Consider the exact sequence of complexes
$$\begin{CD}
\Omega^\bullet_{\X/k}[-1]@>dt\wedge(\cdot )>>\Omega^\bullet_{\X/k}@>>>
\Omega^\bullet_{\X/R}@>>> 0
\end{CD}$$
Applying the exact functor $(\cdot )_{rig}$, we see that
$$\begin{CD}
\Omega^\bullet_{X/k}[-1]@>dt\wedge(\cdot )>>\Omega^\bullet_{X/k}@>>>
\Omega^\bullet_{X/K}@>>> 0
\end{CD}$$ is exact for each value of $\bullet$,
and it is easy to see that these are maps of complexes w.r.t. the
differentials $d$.
\end{proof}

\begin{remark}
The assumption that $X$ is smooth was only used in the application of \cite[7.18-19]{Ni-trace} in the
proof of Lemma \ref{blup}, and can be omitted. In fact, a more natural and elementary way to define the de Rham complex
$(\Omega^\bullet_{X/k},d)$ is the following: if $A$ is an affinoid $K$-algebra, then we define
$d:A\rightarrow \widehat{\Omega}^1_{A/k}$ as the universal continuous $k$-derivation of $A$ into a topological
$A$-module, and
 $\Omega^1_{\mathrm{Sp}\,A/k}$ as the coherent sheaf associated to the finite $A$-module $\widehat{\Omega}^1_{A/k}$ (finiteness
 follows from the fact that, if $A$ is the Tate algebra $K\{x_1,\ldots,x_n\}$, $\widehat{\Omega}^1_{A/k}$ is a
 free $A$-module with basis $dt,dx_1,\ldots,dx_n$ since $k(t)[x_1,\ldots,x_n]$ is dense in $A$).

  One checks that
 this definition behaves well w.r.t. open immersions of $K$-affinoid varieties, so that we can define
 $d:\mathcal{O}_X\rightarrow \Omega^1_{X/k}$ for any rigid $K$-variety $X$. We put
 $\Omega^i_{X/k}=\bigwedge^i\Omega^1_{X/k}$. Then one can prove Lemmas \ref{special}, \ref{etale} and
 \ref{exdef} without smoothness assumptions.
 If $X$ is smooth over $K$ of pure dimension $d$, then $\Omega^1_{X/k}$ is locally
 free of rank $d+1$.
\end{remark}
\subsection{The Gauss-Manin connection}

\begin{definition}\label{def-formcon}
A formal meromorphic connection over $K$ is a $K$-vector space
$M_K$, endowed with a $k$-linear map
$$\partial:M_K\rightarrow M_K$$ satisfying the Leibniz rule
$$\partial(f\cdot m)=\frac{\partial f}{\partial t}\cdot m+f\cdot \partial(m)$$
for each $f$ in $K$ and $m$ in $M_K$, where $\frac{\partial
f}{\partial t}$ is the usual derivative of the Laurent series $f$
w.r.t. the variable $t$.

We say that the formal meromorphic connection
$\partial:M_K\rightarrow M_K$ is regular, if $M_K$ has finite dimension over $K$ and
if there exists an $R$-lattice $M$ in $M_K$ such that $t\partial(M)\subset M$. Such a lattice is
called saturated.
\end{definition}

Let $X$ be a  separated smooth rigid variety over $K$. By
\cite[7.19]{Ni-trace}, we can extend the exact sequence in Lemma
\ref{exdef} to an exact sequence of complexes
$$\begin{CD}
\Omega^\bullet_{X/k}[-2]@>dt\wedge(\cdot )>>\Omega^\bullet_{X/k}[-1]@>dt\wedge(\cdot )>>
\Omega^\bullet_{X/k}@>>> \Omega^\bullet_{X/K}@>>> 0
\end{CD}$$
and we obtain a natural short exact sequence
$$\begin{CD}
0@>>>\Omega^\bullet_{X/K}[-1]@>dt\wedge(\cdot )>>\Omega^\bullet_{X/k}@>>>
\Omega^\bullet_{X/K}@>>> 0
\end{CD}$$
Taking hypercohomology yields a long exact sequence
$$\minCDarrowwidth10pt\begin{CD}
\ldots
@>>>H^{i-1}_{dR}(X/K)@>dt\wedge(\cdot )>>\mathbb{H}^i(X,\Omega^\bullet_{X/k})@>>>
H^{i}_{dR}(X/K)@>\partial_X>> H^{i}_{dR}(X/K)@>>> \ldots
\end{CD}$$
\begin{definition}[Gauss-Manin connection]
We call the connecting homomorphism
$$\partial_X:H^i_{dR}(X/K)\rightarrow H^i_{dR}(X/K)$$ the
Gauss-Manin connection associated to the separated smooth rigid $K$-variety
$X$.
\end{definition}
\begin{lemma}\label{isconn}
The pair $(H^i_{dR}(X/K),\partial_X)$ is a formal meromorphic connection
over $K$, for any separated smooth  rigid $K$-variety $X$ and any integer
$i\geq 0$.
\end{lemma}
\begin{proof}
If $X$ is quasi-Stein, then by Kiehl's Theorem B \cite[2.4]{Kiehl}, $H^i_{dR}(X/K)$ (resp. $H^i_{dR}(X/k)$) is
simply
$\mathcal{H}^i(\Omega^\bullet_{X/K}(X))$ (resp.
$\mathcal{H}^i(\Omega^\bullet_{X/k}(X))$), and the sequence
$$\begin{CD}
0@>>>\Omega^\bullet_{X/K}(X)[-1]@>dt\wedge(\cdot )>>\Omega^\bullet_{X/k}(X)@>>>
\Omega^\bullet_{X/K}(X)@>>> 0
\end{CD}$$
is exact.

Therefore, for any closed $i$-form $\omega$ in
$\Omega^i_{X/K}(X)$, the image of $[\omega]\in H^i_{dR}(X/K)$
under $\partial_X$ can be computed as follows: choose an element
$\omega'$ in $\Omega^i_{X/k}(X)$ mapping to $\omega$. Then
$d\omega'$ maps to zero in $\Omega^{i+1}_{X/K}(X)$, so there
exists an element $\alpha$ in $\Omega^{i}_{X/K}(X)$ with
$d\omega'=dt\wedge \alpha$, and  $\partial_X([\omega])=[\alpha]$.
It is clear that
$$\partial_X(f[\omega])=\frac{\partial f}{\partial
t}[\omega]+f\partial_X([\omega])$$ for each $f$ in $K$.

The general case follows by using \v{C}ech cohomology w.r.t. an admissible
cover by affinoid domains.
\end{proof}

\begin{prop}
If $h:Y\rightarrow X$ is a morphism of separated smooth rigid $K$-varieties,
then the square
$$\begin{CD}
H^{i}_{dR}(X/K)@>\partial_X>> H^{i}_{dR}(X/K)
\\@Vh^*VV @VVh^*V
\\H^{i}_{dR}(Y/K)@>\partial_Y>> H^{i}_{dR}(Y/K)
\end{CD}$$
commutes for each $i\geq 0$.
\end{prop}
\begin{proof}
This is clear from the definition.
\end{proof}


%

\begin{theorem}
If $X$ is a proper and smooth rigid $K$-variety, then for each $i\geq 0$, the pair $(H^i_{dR}(X/K),\partial_X)$
is a regular formal meromorphic connection.
\end{theorem}
The proof will appear in a forthcoming article. It uses a standard argument involving
complexes of logarithmic differential forms.
\section{Comparison with the algebraic setting}\label{subsec-proj2}
Let $S$ be a
smooth $k$-curve, let $0$ be a point in $S(k)$,
and fix a local parameter $t$ on $S$ at $0$. Without loss of generality, we assume that $t$ is defined on $S$.
 We put $S^o=S-\{0\}$. The parameter $t$ defines a morphism
of $k$-schemes $\widehat{\eta}:\Spec K\rightarrow S^o$.

Let $f:Y\rightarrow S^o$ be a smooth and proper morphism, and put $X=Y\times_{S^o} \widehat{\eta}$.
We denote, for each $i\geq 0$,
by $$\nabla_{\!\!\frac{\partial}{\partial t}}:R^if_*(\Omega^\bullet_{Y/S^o})\rightarrow
R^if_*(\Omega^\bullet_{Y/S^o})$$
the covariant derivative w.r.t. $\frac{\partial}{\partial t}$ of
 the algebraic Gauss-Manin connection $\nabla$ associated to $f$ \cite{katz-oda}.
 Recall that, if $k=\C$, the complex analytic
 connection induced by $\nabla$ is the natural connection associated to the local
 subsystem $R^if_*^h(\C_{Y^h})$ of
the locally free $\mathcal{O}_{(S^o)^h}$-module $R^if^h_*(\Omega^\bullet_{Y^h/(S^o)^h})$.

By the base change property for relative de Rham cohomology, we get canonical isomorphisms
$$\widehat{\eta}^*R^if^o_*(\Omega^\bullet_{Y/S^o})\cong H^i_{dR}(X/K)$$
Applying the Leibniz rule, we see that there exists a unique formal meromorphic connection $\partial_t$ on $H^i_{dR}(X/K)$
such that the square
$$\begin{CD}
R^if_*(\Omega^\bullet_{Y/S^o})(U)@>\nabla_{\!\!\frac{\partial}{\partial t}}>>R^if_*(\Omega^\bullet_{Y/S^o})(U)
\\ @VVV @VVV
\\ H_{dR}^i(X/K)@>\partial_t>> H_{dR}^i(X/K)
\end{CD}$$
commutes for every open subscheme $U$ of $S^o$.

\begin{theorem}[Comparison theorem: proper case]\label{compar-prop}
The natural GAGA isomorphism $H^i_{dR}(X/K)\rightarrow H^i_{dR}(X^{an}/K)$ commutes with $\partial_t$
and $\partial_{X^{an}}$.
\end{theorem}
We only give a brief sketch of the proof; details will appear in a forthcoming article. The main difficulty is interpreting the cohomology of the complex
$\Omega^{\bullet}_{X^{an}/k}$
in terms of algebraic geometry. We choose a compactification $g:Z\rightarrow S$ of $f$ over $S$ such that $Z$
is smooth over $k$ and $Z_s=Z\times_S 0$ has strict normal crossings.
Denote by $\mZ$ the $t$-adic completion of $g$. Its
generic fiber $\mZ_\eta$ is canonically isomorphic to $X^{an}$.

Since $(sp_\mZ)_*$ is exact on coherent $\mathcal{O}_{\mZ_\eta}$-modules, Leray's spectral sequence yields isomorphisms
$$\mathbb{H}^i(X^{an},\Omega^{\bullet}_{X^{an}/k})\cong \mathbb{H}^i(\mZ,(sp_{\mZ})_*\Omega^{\bullet}_{X^{an}/k})$$
One defines a logarithmic  subcomplex $\Omega^{\bullet}_{\mZ/k}(\log Z_s)$ of $(sp_\mZ)_*\Omega^{\bullet}_{X^{an}/k}$
in the usual way, and one shows that the inclusion morphism
$$\Omega^{\bullet}_{\mZ/k}(\log Z_s)\rightarrow (sp_\mZ)_*\Omega^{\bullet}_{X^{an}/k}$$
is a quasi-isomorphism. Then one applies the GFGA comparison theorem \cite[4.1.5]{ega3}.
\section{Comparison with the complex analytic setting: the Milnor
fibration} In this section, we suppose that $k=\C$. Let $f:\C^{n+1}\rightarrow \C$ be a complex analytic map, for some integer $n>0$,
 and assume that $f$ has an isolated singularity at $x\in f^{-1}(0)$.
We fix a complex coordinate $t$ on $\C$. Taking the formal completion of $f$
at $x$, we get a generically smooth special formal $R$-scheme
$$\widehat{f}:\X=\mathrm{Spf}\,\widehat{\mathcal{O}}_{\C^{n+1},x}\rightarrow
\mathrm{Spf}\,R$$ We call its generic fiber $\X_\eta$ the analytic Milnor fiber of $f$ at $x$ (see \cite{NiSe-Milnor}
\cite{NiSe}),
and denote it by $\mathscr{F}_x$. It is a complete invariant of the formal
germ of $f$ at $x$ \cite[8.5]{Ni-trace}. Note that, since $\X$ is affine, we have
$$R^i\widehat{f}_*(\Omega^\bullet_{\X/R})=\mathcal{H}^i(\Omega^\bullet_{\X/R}(\X))$$ for each $i\geq 0$.

Let $B=B(x,\varepsilon)\subset \C^{n+1}$
be the open ball around $x$ with radius $\varepsilon$, and let $D=B(0,\eta)\subset \C$ be the open disc
around $0$ with radius $\eta$. We denote by $D^o$ the punctured disc $D-\{0\}$.
If $0<\eta\ll \varepsilon \ll 1$, then
the induced map $$f^o:X^o=B\cap f^{-1}(D^o)\rightarrow D^o$$ is a locally
trivial fibration;
this is the so-called Milnor fibration of $f$ at $x$. We put $X'=B\cap f^{-1}(D)$ and denote by $f':X'\rightarrow D$
the restriction of $f$.

We consider the universal covering
 $$\widetilde{D^o}=\{z\in \C\,|\,\Im(z)>-\frac{\log \eta}{2\pi}\}\rightarrow D^o:z\mapsto \exp(2\pi i z)$$
 The universal fiber
$F_x:=X^o\times_{D^o}\widetilde{D^o}$ of $f^o$ is called the
 Milnor fiber of $f$ at $x$. The group $\pi_1(D^o)=\Z$ of covering transformations acts on the singular cohomology spaces
 $H^i(F_x,\Z)$. The action of the canonical generator $z\mapsto z+1$ is called the monodromy transformation.

 For each $i\geq
0$, the sheaf $R^i(f^o)_*(\C_{X^o})$ is a local system on $D^o$. It corresponds to
the monodromy representation $\pi_1(D^o)\rightarrow Aut_{\C}(H^i(F_x,\C))$.
It vanishes for $i\notin \{0,n\}$, and its fiber is a complex
vector space of dimension $\mu$ (resp. $1$) for $i=n$ (resp.
$i=0$), where $\mu$ is the Milnor number of $f$ at $x$. Since
$$R^i(f^o)_*(\Omega^\bullet_{X^{o}/D^{o}})\cong R^i(f^o)_*(\C_{X^o})\otimes_{\C_{D^o}}\mathcal{O}_{D^{o}}$$
for $i\geq 0$, the local system $R^i(f^o)_*(\C_{X^o})$ defines an
integrable connection $\nabla$ on the relative de Rham
cohomology sheaf $R^i(f^o)_*(\Omega^\bullet_{X^{o}/D^{o}})$,
called the Gauss-Manin connection \cite[2.29]{deligne}.

Brieskorn \cite[Prop.\,1.6+Satz\,1]{brieskorn}  proved that $R^if'_*(\Omega^\bullet_{X'/D})$ is a coherent
$\mathcal{O}_D$-module for $i\geq 0$, whose stalk
$H^i_{f,x}:=(R^if'_*(\Omega^\bullet_{X'/D}))_0$ at $0$ is canonically isomorphic
to the $i$-th cohomology space $\mathcal{H}^i(\Omega^{\bullet}_{X'/D,x})$ of the localized relative
de Rham complex of $f$ at $x$.
Sebastiani
\cite[Cor.\,1]{sebastiani} proved that $H^n_{f,x}$ has no
$t$-torsion; hence, it is a free $\mathcal{O}_{D,0}$-module of
rank $\mu$. As explained in \cite[3.5.1]{Kuli}, the presence of a
connection on $H^i_{f,x}$ for $0\leq i<n$ implies that it is a
free $\mathcal{O}_{D,0}$-module, so it vanishes for $0<i<n$.
We denote by $\widehat{H}^i_{f,x}$ the $t$-adic completion of
$H^i_{f,x}$.
The following proposition is proven as in \cite[3.2]{brieskorn}.
\begin{prop}
For each $i\geq 0$, there exists a canonical isomorphism of
$R$-modules
$$\widehat{H}^i_{f,x}\cong R^i\widehat{f}_*(\Omega^{\bullet}_{\X/R}) $$
In particular, the right hand side is a free $R$-module of rank
$\mu$ (resp. $1$) for $i=n$ (resp. $i=0$), and vanishes for
$i\notin\{0,n\}$.
\end{prop}
All proofs of the coherence of $\widehat{H}^i_{f,x}$ known to me use
either an embedding in a proper family or the theory of
nuclear Fr\'echet algebras (Kiehl-Verdier Theorem); it would be
more satisfactory to find a purely algebraic proof in terms of the
$R$-algebra $\widehat{\mathcal{O}}_{X,x}$.


We can endow $\widehat{H}^i_{f,x}\otimes_R K$ with a formal meromorphic connection $\partial_t$ as follows.
For $i=0$, $\widehat{H}^0_{f,x}$ is canonically isomorphic to $R$, and we endow
$\widehat{H}^0_{f,x}\otimes_R K=K$ with the trivial connection $\partial_t(f)=\partial f/\partial t$.
Now let us consider the remaing case $i=n$.
Since $\Omega^{n+1}_{\mathscr{F}_x/K}=0$, the results in Section \ref{sec-gm} show that the wedge product $dt\wedge
(\cdot )$ induces an isomorphism of coherent
$\mathcal{O}_{\mathscr{F}_x}$-modules
$$dt\wedge:\Omega^{n}_{\mathscr{F}_x/K}\rightarrow
\Omega^{n+1}_{\mathscr{F}_x/\C}=(\Omega^{n+1}_{\X/\C})_{rig}$$ The inverse image of an element
$\omega\in\Omega^{n+1}_{\mathscr{F}_x/\C}(\mathscr{F}_x)$ w.r.t. the
isomorphism $dt\wedge$ is denoted by $\omega/dt\in
\Omega^{n}_{\mathscr{F}_x/K}(\mathscr{F}_x)$. By analogy with the
complex analytic setting, we call $\omega/dt$ the Gelfand-Leray
form associated to $\omega$. We've noted in \cite[2.8]{Ni-trace} that
the natural maps $\Omega^n_{\X/R}(\X)\otimes_R K\rightarrow
\Omega^{n}_{\mathscr{F}_x/K}(\mathscr{F}_x)$ and
$\Omega^{n+1}_{\X/\C}(\X)\otimes_R K\rightarrow
\Omega^{n+1}_{\mathscr{F}_x/\C}(\mathscr{F}_x)$ are injections, and
by \cite[7.20]{Ni-trace}, $ dt\wedge$ restricts to an isomorphism of
$(\widehat{\mathcal{O}}_{X,x}\otimes_R K)$-modules
$$ dt \wedge:\Omega^{n}_{\X/R}(\X)\otimes_R K\rightarrow
\Omega^{n+1}_{\X/\C}(\X)\otimes_R K$$ Now let $a$ be an element of
$\widehat{H}^n_{f,x}$, and choose a representant $\alpha$ of $a$
in $\Omega^{n}_{\X/\C}(\X)$. It is easily seen that the class of
$d\alpha/dt$ in $\widehat{H}^n_{f,x}\otimes_R K$ only depends on
$a$. The map
$$\partial_t:\widehat{H}^n_{f,x}\rightarrow
\widehat{H}^n_{f,x}\otimes_R K:a\mapsto d\alpha/dt$$ satisfies the
Leibniz rule, and hence extends uniquely to a formal meromorphic
connection $\partial_t$ on
$\widehat{H}^n_{f,x}\otimes_R K$.

We'll now show how $\partial_t$ relates to the Gauss-Manin connection
$\nabla$ on $R^if^o_*(\Omega^\bullet_{X^o/S^o})$. We put $\mathscr{R}=\mathcal{O}_{\C,0}$ and we
denote by $\mathscr{K}$ its quotient field. If $V$ is a vector space over $\mathscr{K}$, the notions of (regular)
meromorphic connection and saturated lattice
are defined similarily as in Definition \ref{def-formcon}. We put
$$\mathbb{K}=\lim_{\stackrel{\longrightarrow}{E}}\mathcal{O}_{\C}(E-\{0\})$$ where $E$ runs through
a fundamental system of open neighbourhoods of $0$ in $\C$. For each $i\geq 0$, there is a canonical isomorphism
$$H^i_{f,x}\otimes_{\mathscr{R}}\mathbb{K}\cong\lim_{\stackrel{\longrightarrow}{E}}R^if^o_*(\Omega^\bullet_{X^o/S^o})(E-\{0\})$$
and by passing to the limit,
the covariant derivative $\nabla_{\!\!\frac{\partial}{\partial t}}$ defines a $\mathbb{C}$-linear
endomorphism $\overline{\partial}_t$ of $H^i_{f,x}\otimes_{\mathscr{R}}\mathbb{K}$.
\begin{theorem}
For $i=0$ and $i=n$, there exists a unique meromorphic connection $\widetilde{\partial}_t$ on
$H_{f,x}^i\otimes_{\mathscr{R}}\mathscr{K}$
such that the natural maps of $\mathscr{K}$-vector spaces
\begin{eqnarray*}
H_{f,x}^i\otimes_{\mathscr{R}}\mathscr{K}&\rightarrow &\widehat{H}^i_{f,x}\otimes_R K
\\H_{f,x}^i\otimes_{\mathscr{R}}\mathscr{K}&\rightarrow &H^i_{f,x}\otimes_{\mathscr{R}}\mathbb{K}
\end{eqnarray*}
commute with $\widetilde{\partial}_t$, $\partial_t$ and $\overline{\partial}_t$. The connections $\widetilde{\partial}_t$
and $\partial_t$ are regular.
\end{theorem}
\begin{proof}
Uniqueness is clear.
 Existence for $i=0$ follows from the isomorphism $R^0f^o_*(\C_{X^o})\cong\C_{D^o}$: we have $H^i_{f,x}=\mathscr{R}$,
 and $\widetilde{\partial}_t$ is the trivial meromorphic connection on $\mathscr{K}$.
The case $i=n$ follows from Brieskorn's result
\cite[Satz\,1]{brieskorn}: the meromorphic connection $\widetilde{\partial}_t$ is the one constructed there.
It is well-known that $\widetilde{\partial}_t$ is regular \cite[I.8.3]{Kuli}. If $M$ is a saturated lattice
for $\widetilde{\partial}_t$, then $M\otimes_{\mathscr{R}}R$ is a saturated lattice for $\partial_t$, so $\partial_t$
is regular as well.
\end{proof}

In the light of the motivic monodromy
conjecture \cite[\S 2.4]{DL5}, it is highly intriguing that
the local motivic zeta function of $f$ at $x$ can be expressed in
terms of certain motivic integrals of the Gelfand-Leray form
$\omega/dt$, where $\omega$ is an arbitrary gauge form on $\X$
over $\C$ (see \cite[9.7]{Ni-trace}).

 It follows from Berthelot's construction \cite[0.2.6]{bert} that
$\mathscr{F}_x$ is a partially proper quasi-Stein space, so
$H_{dR}^i(\mathscr{F}_x/K)$ coincides with the overconvergent de Rham cohomology from \cite{rigiddeRham}
and is simply given by
$\mathcal{H}^i(\Omega^{\bullet}_{\mathscr{F}_x/K}(\mathscr{F}_x))$, by
\cite[1.8(b)]{rigiddeRham} and Kiehl's Theorem B
\cite[2.4]{Kiehl}. The arguments in
\cite{finitedeRham} can be used to show that $H_{dR}^i(\mathscr{F}_x/K)$ has finite dimension
for each $i\geq 0$. The natural map of complexes
$$\Omega^\bullet_{\X/R}(\X)\rightarrow \Omega^\bullet_{\mathscr{F}_x/K}(\mathscr{F}_x)$$
induces a natural map of $K$-vector spaces
$$\gamma_i:\widehat{H}^i_{f,x}\otimes_R K\rightarrow H^i_{dR}(\mathscr{F}_x/K) $$ for each $i\geq 0$.

\begin{prop}
The natural map
$$\gamma_0:\widehat{H}^0_{f,x}\otimes_R K\rightarrow H^0_{dR}(\mathscr{F}_x/K) $$ is an isomorphism, and
$H^0_{dR}(\mathscr{F}_x/K)=K$.
 Moreover, the natural diagram
$$\begin{CD}
\widehat{H}^i_{f,x}\otimes_R K@>\gamma_i>> H^i_{dR}(\mathscr{F}_x/K)
\\ @V\partial_t VV @V\partial_{\mathscr{F}_x}VV
\\ \widehat{H}^i_{f,x}\otimes_R K@>\gamma_i>> H^i_{dR}(\mathscr{F}_x/K)
\end{CD}$$
 commutes for $i=0$ and $i=n$.
\end{prop}
\begin{proof}
To show that $\gamma_0$ is an isomorphism, it suffices to show that $\mathscr{F}_x$ is geometrically connected.
This follows from the equality of the $\ell$-adic Betti numbers of $\mathscr{F}_x$ and the Betti numbers of $F_x$ \cite[9.2]{NiSe}.
 Since $\mathscr{F}_x$ is quasi-Stein, the
proof of Lemma \ref{isconn} shows that the above diagram commutes for $i=0,n$.
\end{proof}
\begin{conjecture}\label{conj1}
The natural map
$$\widehat{H}^i_{f,x}\otimes_R K\rightarrow H^i_{dR}(\mathscr{F}_x/K) $$
is an isomorphism for each $i\geq 0$. In particular, the right
hand side vanishes unless $i=0$ or $i=n$, and
 has dimension $\mu$ (resp. $1$) for $i=n$ (resp. $i=0$).
\end{conjecture}
This conjecture roughly states that the de Rham cohomology of
$\mathscr{F}_x$ can be computed using only differential forms on
$\mathscr{F}_x$ which are meromorphic along the special fiber of
$\X$.
Let us give some motivation for the conjecture. Heuristically, it fits with
 similar results in complex GAGA
(computation of de Rham cohomology using meromorphic or
logarithmic forms), which also appear in the proof of Theorem \ref{compar-prop} (note, however, that
the isomorphism $$R^ig_*(\Omega^\bullet_{\mY/R})\otimes_R K\cong H^i_{dR}(\mY_\eta/K)$$ is easy to prove when
$g:\mY\rightarrow \Spf R$ is $stft$, using exactness of $(sp_{\mY})_*$ on coherent $\mathcal{O}_{\mY_\eta}$-modules
 and the fact that $(sp_{\mY})_*\Omega^{\bullet}_{\mY_\eta/K}\cong \Omega^\bullet_{\mY/R}\otimes_R K$).

 Moreover, as we showed in \cite[9.2]{NiSe}, the
$i$-th \'etale $\ell$-adic Betti number of $\mathscr{F}_x$ vanishes
for $i\notin \{0,n\}$ and equals $\mu$ (resp. $1$) for
$i=n$ (resp. $i=0$). More precisely, there exists for each $i\geq 0$ and each prime $\ell$ an isomorphism
$$H^i(\mathscr{F}_x\widehat{\times}_K \widehat{K^a},\Q_\ell)\cong H^i(F_x,\Q_\ell)$$ which identifies
the action of the canonical generator of $G(K^a/K)$ on the left hand side, with the monodromy transformation on the
right hand side (here $\widehat{K^a}$ denotes the completion of an algebraic closure of $K$). Finally, we proved in
\cite{Ni-sing} that the rational singular cohomology of $\mathscr{F}_x\widehat{\times}_K \widehat{K^a}$ (viewed as
a Berkovich space) is isomorphic to the weight zero part of Steenbrink's \cite{steenbrink-vanish} mixed Hodge structure on the rational singular
cohomology of
$F_x$. In particular, it is concentrated in degree $0$ and $n$.


\bibliographystyle{hplain}
\bibliography{wanbib,wanbib2}
\end{document}